\documentclass[11pt,reqno]{amsart}
\usepackage[margin=1in]{geometry}
\usepackage{amsmath,amssymb,amsthm}
\usepackage{enumerate}
\usepackage{cite}
\usepackage[colorlinks=true,linkcolor=blue,citecolor=blue,urlcolor=blue]{hyperref}

\newtheorem{theorem}{Theorem}[section]
\newtheorem{lemma}[theorem]{Lemma}
\newtheorem{proposition}{Proposition}
\theoremstyle{definition}
\newtheorem{example}{Example}
\theoremstyle{remark}

\numberwithin{equation}{section}

\newcommand{\omex}{\overline{\mathrm{mex}}}          
\newcommand{\sm}[1]{\sigma_{#1}^{\overline{\mathrm{mex}}}}
\newcommand{\cP}{\overline{\mathcal{P}}}
\newcommand{\op}{\overline{p}}

\title{The minimal excludant of $\ell$-regular overpartitions}
\author{S. Bhowmick}
\address{Sourav Bhowmick, Department of Mathematics, National Institute of Technology, Raipur, Chhattisgarh 492010.}
\email{souravbhowmick578@gmail.com, sbhowmick.phd2025.maths@nitrr.ac.in}

\author{N. K. Meher}
\address{Nabin Kumar Meher, Department of Mathematics, National Institute of Technology, Raipur, Chhattisgarh 492010.}
\email{mehernabin@gmail.com, nkmeher.maths@nitrr.ac.in}

\begin{document}

\begin{abstract}
	Andrews and Newman \cite{AN1} introduced the minimal excludant (mex) of an integer partition, and Aricheta and Donato \cite{AD} extended it to overpartitions. More recently, Meher \cite{Meher} and Tejash, Bhowmick and Meher \cite{TBM} studied the analogous statistic for $\ell$-regular partitions. In this article, for an integer $\ell\ge 2$, we define the minimal excludant $\omex_\ell(\pi)$ of an $\ell$-regular overpartition $\pi$ as the smallest positive integer that is not divisible by $\ell$ and does not occur as a non-overlined part of $\pi$. We study $\sm{\ell}(n)$, the sum of $\omex_\ell(\pi)$ over all $\ell$-regular overpartitions $\pi$ of $n$. We obtain generating functions for $\sm{\ell}(n)$ as  eta-quotients for $\ell=2$ and $\ell=3$. Later, motivated from \cite{TBM}, we obtained explicit generating functions for $\sm{\ell}(n)$.  We show that $\sm{\ell}(n)$ is weakly increasing in $n$, and we derive Hardy--Ramanujan type asymptotic formulas for $\sm{2}(n)$ and $\sm{3}(n)$. We then determine the parity of both functions completely: $\sm{2}(n)$ is even for every $n\ge 1$, while $\sm{3}(n)$ is odd if and only if $n=0$ or $n$ is a centred hexagonal number. Also, we obtain for every $\ell\ge3$, $\sm{\ell}(n)$ is odd for infinitely many $n$, and the number of such $n\le X$ is $\gg\sqrt{X}$. Finally, we establish congruences for these functions modulo $3$, $4$ and $8$.
\end{abstract}

\maketitle

\section{Introduction and Statement of Results}\label{sec:intro}

A partition $\lambda=(\lambda_1,\lambda_2,\ldots,\lambda_k)$ of a non-negative integer $n$ is a non-increasing sequence of positive integers summing to $n$. Let $p(n)$ denote the number of partitions of $n$. Throughout, for $|q|<1$ we write $(a;q)_\infty:=\prod_{k\ge0}(1-aq^k)$, and for brevity we set $f_r:=(q^r;q^r)_\infty$ for every positive integer $r$. The generating function of $p(n)$ is given by
\[
\sum_{n\ge0}p(n)q^n=\prod_{r\ge1}\frac{1}{1-q^r}=\frac{1}{f_1}.
\]
The arithmetic and analytic properties of $p(n)$ have been studied extensively; see, for example, \cite{HR,Ram,GO,Ono,AO}. Grabner and Knopfmacher \cite{GK} introduced and studied the \emph{smallest gap} of a partition, namely the least positive integer missing from its parts. In 2019, Andrews and Newman \cite{AN1,AN2} rediscovered this statistic under the name minimal excludant, $\mathrm{mex}(\lambda)$. They defined $\sigma\mathrm{mex}(n):=\sum_{\lambda\vdash n}\mathrm{mex}(\lambda)$ and proved
\begin{equation}\label{eq:1.1}
\sigma\mathrm{mex}(n)=D_2(n),
\end{equation}
where $D_2(n)$ is the number of partitions of $n$ into distinct parts using two colours. Ballantine and Merca \cite{BM} gave a combinatorial proof of \eqref{eq:1.1}. Barman and Singh \cite{BS1,BS2} and da Silva and Sellers \cite{dSS} established congruences and parity results for mex-related partition functions.

An \emph{overpartition} of $n$ is a partition of $n$ in which the first occurrence of each distinct part may be overlined. Let $\op(n)$ be the number of overpartitions of $n$. For example, $\op(3)=8$:
\[
3,\ \overline{3},\ 2+1,\ \overline{2}+1,\ 2+\overline{1},\ \overline{2}+\overline{1},\ 1+1+1,\ \overline{1}+1+1.
\]
The generating function of $\op(n)$ is given by
\[
\sum_{n\ge0}\op(n)q^n=\prod_{r\ge1}\frac{1+q^r}{1-q^r}=\frac{(-q;q)_\infty}{(q;q)_\infty}=\frac{f_2}{f_1^2}.
\]
Overpartitions were studied systematically by Corteel and Lovejoy \cite{CL}, and the arithmetic properties of $\op(n)$ have been investigated by many authors; see, for instance, \cite{Mah,HS05}. Aricheta and Donato \cite{AD} defined the minimal excludant of an overpartition $\pi$ as the smallest positive integer that is not a part of the non-overlined parts of $\pi$. They proved that $\overline{\sigma}\mathrm{mex}(n)=D_3(n)$, the number of partitions of $n$ into distinct parts using three colours. They also showed that $\overline{\sigma}\mathrm{mex}(n)$ is odd exactly when $n$ is a triangular number.

A partition is called \emph{$\ell$-regular} if none of its parts is divisible by $\ell$. Let $p_\ell(n)$ denote the number of $\ell$-regular partitions of $n$. Then the generating function of $p_\ell(n)$ is given by
\[
\sum_{n\ge0}p_\ell(n)q^n=\prod_{\substack{r\ge1\\ \ell\nmid r}}\frac{1}{1-q^r}=\frac{f_\ell}{f_1}.
\]
Divisibility, congruence and parity properties of $p_\ell(n)$ have been studied, for example, in \cite{GO,LP,Pen,HS10}. Recently Meher \cite{Meher} defined the minimal excludant of an $\ell$-regular partition $\lambda$ by $\mathrm{mex}_\ell(\lambda):=\min\{r\ge1:\ \ell\nmid r,\ r\notin\lambda\}$. He found the generating functions of $\sigma_\ell^{\mathrm{mex}}(n)$ for $\ell=2,3$, together with asymptotic formulas. He showed that $\sigma_2^{\mathrm{mex}}(n)$ is odd exactly at the generalized pentagonal numbers, and that $\sigma_3^{\mathrm{mex}}(n)$ is odd for $\gg\log\log X$ values of $n\le X$. Motivated by these results, we study the minimal excludant of $\ell$-regular overpartitions.

Let $\ell\ge2$ be a positive integer. An overpartition is \emph{$\ell$-regular} if none of its parts, overlined or not, is divisible by $\ell$. Let $\cP_\ell(n)$ denote the set of all $\ell$-regular overpartitions of $n$, and $\op_\ell(n):=|\cP_\ell(n)|$. The overlined parts of an $\ell$-regular overpartition form a partition into distinct parts not divisible by $\ell$, and the non-overlined parts form an $\ell$-regular partition. Since
\[
\prod_{\substack{r\ge1\\ \ell\nmid r}}(1+q^r)=\frac{(-q;q)_\infty}{(-q^\ell;q^\ell)_\infty}=\frac{f_2}{f_1}\cdot\frac{f_\ell}{f_{2\ell}},
\]
the generating function of $\op_\ell(n)$ is given by
\[
\overline{P}_\ell(q):=\sum_{n\ge0}\op_\ell(n)q^n=\prod_{\substack{r\ge1\\ \ell\nmid r}}\frac{1+q^r}{1-q^r}=\frac{f_2f_\ell^2}{f_1^2f_{2\ell}}.
\]
In particular, $\overline{P}_2(q)=\dfrac{f_2^3}{f_1^2f_4}$ and $\overline{P}_3(q)=\dfrac{f_2f_3^2}{f_1^2f_6}$. The function $\op_\ell(n)$ was introduced by Lovejoy \cite{Lov03}, and its arithmetic properties have been studied in \cite{Shen,AMS,Chern,HaoShen}.
For $\pi\in\cP_\ell(n)$, the \emph{$\ell$-regular minimal excludant} is defined by
\[
\omex_\ell(\pi):=\min\{r\ge1:\ \ell\nmid r,\ r \text{ is not a non-overlined part of } \pi\}.
\]
In other words, $\omex_\ell(\pi)$ is the smallest positive integer that is not divisible by $\ell$ and does not occur as a non-overlined part of $\pi$. Overlined parts play no role in $\omex_\ell(\pi)$. We call a positive integer \emph{admissible} if it is not divisible by $\ell$. Thus every part of an $\ell$-regular overpartition, overlined or not, is admissible, and $\omex_\ell(\pi)$ is the smallest admissible integer that does not occur as a non-overlined part of $\pi$. Let $m$ be a positive integer such that $\ell\nmid m$. We define
\[
p_\ell^{\overline{\mathrm{mex}}}(m,n):=\#\{\pi\in\cP_\ell(n):\ \omex_\ell(\pi)=m\}.
\]
If $\ell\mid m$, then $p_\ell^{\overline{\mathrm{mex}}}(m,n)=0$ by convention. We consider some examples to understand the above definitions.

\begin{example}
Let $\ell=4$ and $n=6$. The admissible parts are $1,2,3,5,6,7,9,10,11,\ldots$.
\begin{enumerate}[(1)]
\item Let $\pi=6$. The only non-overlined part is $6$, so $1$ does not occur as a non-overlined part. Hence $\omex_4(\pi)=1$.
\item Let $\pi=5+1$. The non-overlined parts contain $1$ but not $2$, so $\omex_4(\pi)=2$. On the other hand, for $\pi=5+\overline{1}$ the part $1$ is overlined and hence ignored, so $\omex_4(\pi)=1$.
\item Let $\pi=2+2+1+1$. The non-overlined parts contain $1,2$ but not $3$, so $\omex_4(\pi)=3$.
\item Let $\pi=3+2+1$. The non-overlined parts contain $1,2,3$. The integer $4$ is not admissible, and $5$ does not occur, hence $\omex_4(\pi)=5$. For $\pi=\overline{3}+2+1$, the non-overlined parts are $2,1$, so $\omex_4(\pi)=3$.
\item Let $\pi=3+\overline{2}+1$. The non-overlined parts are $3,1$, so $\omex_4(\pi)=2$.
\end{enumerate}
\end{example}

Our main object of study is the function
\[
\sm{\ell}(n):=\sum_{\pi\in\cP_\ell(n)}\omex_\ell(\pi),\qquad \sm{\ell}(0):=1,
\]
the natural analogue of $\sigma_\ell^{\mathrm{mex}}(n)$ for $\ell$-regular overpartitions. We obtain explicit generating functions for $\sm{\ell}(n)$, expressed as eta-quotients, for $\ell=2$ and $\ell=3$.

\begin{theorem}\label{thm:1.1}
We have
\begin{equation}\label{eq:1.2}
\sum_{n\ge0}\sm{2}(n)q^n=\frac{(q^2;q^2)_\infty^8}{(q;q)_\infty^4(q^4;q^4)_\infty^3}=\frac{f_2^8}{f_1^4f_4^3}.
\end{equation}
\end{theorem}

\begin{theorem}\label{thm:1.2}
We have
\begin{equation}\label{eq:1.3}
\sum_{n\ge0}\sm{3}(n)q^n=\frac{f_2f_6^4}{f_1^2f_{12}^2}+q\,\frac{f_2f_3^2f_{12}^2}{f_1^2f_6^2}.
\end{equation}
\end{theorem}

In general, motivated by \cite{TBM} we derive the following.

\begin{theorem}\label{thm:1.3}
Let $\ell\ge2$ be an integer. If $\ell$ is even, then
\begin{equation}\label{eq:1.4}
\sum_{n\ge0}\sm{\ell}(n)q^n=\frac{f_2f_\ell^2f_{\ell(\ell-1)}}{f_1^2f_{2\ell}}\sum_{s=1}^{\ell/2}q^{\frac{s(s-1)}{2}}\Big(-q^{\frac{\ell(\ell+2s-3)}{2}};q^{\ell(\ell-1)}\Big)_\infty\Big(-q^{\frac{\ell(\ell-2s+1)}{2}};q^{\ell(\ell-1)}\Big)_\infty.
\end{equation}
If $\ell$ is odd, then
\begin{multline}\label{eq:1.5}
\sum_{n\ge0}\sm{\ell}(n)q^n=\frac{f_2f_\ell^2f_{\ell(\ell-1)}}{f_1^2f_{2\ell}}\Bigg[\sum_{s=1}^{(\ell-1)/2}q^{\frac{s(s-1)}{2}}\Big(-q^{\frac{\ell(\ell+2s-3)}{2}};q^{\ell(\ell-1)}\Big)_\infty\Big(-q^{\frac{\ell(\ell-2s+1)}{2}};q^{\ell(\ell-1)}\Big)_\infty\\
+q^{\frac{\ell^2-1}{8}}\,\frac{f_{2\ell(\ell-1)}^2}{f_{\ell(\ell-1)}^2}\Bigg].
\end{multline}
\end{theorem}

Further, we show that $\{\sm{\ell}(n)\}_{n\ge0}$ is weakly increasing for every $\ell\ge2$. In particular, we have the following result.

\begin{theorem}\label{thm:1.4}
The sequence $\{\sm{\ell}(n)\}_{n\ge0}$ is weakly increasing. Moreover, it satisfies the inequality
\begin{equation}\label{eq:1.6}
\sm{\ell}(n+1)-\sm{\ell}(n)\ge\op_\ell(n)-\op_\ell(n-1).
\end{equation}
\end{theorem}

Using Ingham's Tauberian theorem, we establish Hardy--Ramanujan type asymptotic formulas for $\sm{2}(n)$ and $\sm{3}(n)$.

\begin{theorem}\label{thm:1.5}
We have
\[
\sm{2}(n)\sim\frac{1}{2\sqrt{2n}}\,e^{\pi\sqrt{n/2}},\qquad\text{as } n\to\infty.
\]
\end{theorem}

\begin{theorem}\label{thm:1.6}
We have
\[
\sm{3}(n)\sim\frac{1}{4\sqrt{n}}\,e^{\pi\sqrt{2n/3}},\qquad\text{as } n\to\infty.
\]
\end{theorem}

We further determine the parity of $\sm{2}(n)$ and $\sm{3}(n)$ completely, and deduce the number of $n\le X$ for which $\sm{3}(n)$ is odd and even. In particular, we get the following results.

\begin{theorem}\label{thm:1.7}
\leavevmode
\begin{enumerate}[(1)]
\item $\sm{2}(n)$ is an even integer for every $n\ge1$.
\item $\sm{3}(n)\equiv1\pmod 2$ if and only if $n=0$ or $n=3k^2+3k+1$ for some non-negative integer $k$.
\end{enumerate}
\end{theorem}

\begin{theorem}\label{thm:1.8}
For every real number $X\ge1$, we have
\[
\#\{n\le X:\ \sm{3}(n)\text{ is an odd integer}\}=\sqrt{\frac{X}{3}}+O(1),
\]
while
\[
\#\{n\le X:\ \sm{3}(n)\text{ is an even integer}\}=X-\sqrt{\frac{X}{3}}+O(1).
\]
\end{theorem}

More generally, for every $\ell\ge2$ we determine the parity of $\sm{\ell}(n)$ completely.

\begin{theorem}\label{thm:parity-general}
Let $\ell\ge2$. Then $\sm{\ell}(n)$ is odd if and only if $n=0$ or
\[
n=R_\ell(\ell k+s)=\frac{\ell(\ell-1)}{2}k^2+(s-1)\ell k+\frac{s(s-1)}{2}
\]
for some integers $k\ge0$ and $2\le s\le\ell-1$. Consequently, for every real number $X\ge1$,
\[
\#\{0\le n\le X:\ \sm{\ell}(n)\ \text{is odd}\}=(\ell-2)\sqrt{\frac{2X}{\ell(\ell-1)}}+O(1).
\]
In particular, for every $\ell\ge3$, $\sm{\ell}(n)$ is odd for infinitely many $n$, and the number of such $n\le X$ is $\gg\sqrt{X}$.
\end{theorem}

Finally, we prove the following congruences.

\begin{theorem}\label{thm:1.9}
For every non-negative integer $n$,
\begin{enumerate}[(1)]
\item $\sm{2}(3n+2)\equiv0\pmod3$;
\item for $n\ge1$, $\sm{2}(n)\equiv2\pmod4$ if $n=2k^2$ for some $k\ge1$, and $\sm{2}(n)\equiv0\pmod4$;
\item $\displaystyle\sum_{n\ge0}\sm{2}(2n+1)q^n=4\frac{f_4^4}{f_1^2f_2}$. Consequently, for odd $m$, $\sm{2}(m)\equiv4\pmod8$ if $m$ is a perfect square, and $\sm{2}(m)\equiv0\pmod8$;
\item for $n\ge1$, $\sm{3}(3n)\equiv2\pmod4$ if $n=3k^2$ for some $k\ge1$, and $\sm{3}(3n)\equiv0\pmod4$ otherwise. In particular $\sm{3}(9n+3)\equiv\sm{3}(9n+6)\equiv0\pmod4$.
\end{enumerate}
\end{theorem}

\section{Preliminaries}

\begin{lemma}\label{lem:2.1}
We have
\[
\sum_{t\in\mathbb{Z}}q^{t^2}=\frac{(q^2;q^2)_\infty^5}{(q;q)_\infty^2(q^4;q^4)_\infty^2}=\frac{f_2^5}{f_1^2f_4^2}.
\]
\end{lemma}
\begin{proof}
See \cite[Eq.~(1.5.1)]{Hir}.
\end{proof}

\begin{lemma}\label{lem:2.2}
For $z\ne0$ and $|Q|<1$, we have
\begin{equation}\label{eq:2.1}
\sum_{r=-\infty}^{\infty}Q^{r^2}z^r=\big(-Qz;Q^2\big)_\infty\Big(-\frac{Q}{z};Q^2\Big)_\infty\big(Q^2;Q^2\big)_\infty.
\end{equation}
\end{lemma}
\begin{proof}
See \cite[Theorem 2.8]{And}.
\end{proof}

\begin{lemma}[\cite{Ingham}]\label{lem:2.3}
Let $B(q):=\sum_{n=0}^\infty b(n)q^n$ be a power series with radius of convergence $1$. Assume that $\{b(n)\}$ is a weakly increasing sequence of non-negative real numbers. If there are constants $\beta,\gamma\in\mathbb{R}$ and $\mu>0$ such that
\[
B(e^{-y})\sim\beta y^\gamma e^{\mu/y},\qquad\text{as } y\to0^+,
\]
then
\[
b(n)\sim\frac{\beta\,\mu^{\frac{2\gamma+1}{4}}}{2\sqrt{\pi}\,n^{\frac{2\gamma+3}{4}}}\,e^{2\sqrt{\mu n}},\qquad\text{as } n\to\infty.
\]
\end{lemma}

\begin{lemma}\label{lem:2.4}
For any positive integers $u$ and $j$, and a prime number $p$, we have
\begin{equation}\label{eq:2.2}
(q^u;q^u)_\infty^{p^j}\equiv(q^{pu};q^{pu})_\infty^{p^{j-1}}\pmod{p^j}.
\end{equation}
\end{lemma}

\begin{lemma}\label{lem:2.6}
We have
\begin{equation}\label{eq:2.4}
f_1^3=(q;q)_\infty^3=\sum_{n=0}^\infty(-1)^n(2n+1)q^{\frac{n(n+1)}{2}}.
\end{equation}
\end{lemma}
\begin{proof}
See \cite[Eq.~(1.7.1)]{Hir}.
\end{proof}

\begin{lemma}\label{lem:2.8}
Let $k\ge0$ and $1\le s\le\ell-1$. For a positive integer $\ell\ge2$, let
\[
R_\ell(m)=\sum_{\substack{1\le j<m\\ \ell\nmid j}}j.
\]
Then
\[
R_\ell(\ell k+s)=\frac{\ell(\ell-1)}{2}k^2+(s-1)\ell k+\frac{s(s-1)}{2}.
\]
\end{lemma}
\begin{proof}
See \cite[Lemma 2.4]{TBM}.
\end{proof}

\begin{lemma}\label{lem:2.9}
Let $\ell\ge2$ be an integer. Then we have two cases. If $\ell$ is even, then
\[
\sum_{n\ge0}\sigma_\ell^{\mathrm{mex}}(n)q^n=\frac{f_\ell\,f_{\ell(\ell-1)}}{f_1}\sum_{s=1}^{\ell/2}q^{\frac{s(s-1)}{2}}\Big(-q^{\frac{\ell(\ell+2s-3)}{2}};q^{\ell(\ell-1)}\Big)_\infty\Big(-q^{\frac{\ell(\ell-2s+1)}{2}};q^{\ell(\ell-1)}\Big)_\infty.
\]
If $\ell$ is odd, then
\begin{multline*}
\sum_{n\ge0}\sigma_\ell^{\mathrm{mex}}(n)q^n=\frac{f_\ell\,f_{\ell(\ell-1)}}{f_1}\Bigg[\sum_{s=1}^{(\ell-1)/2}q^{\frac{s(s-1)}{2}}\Big(-q^{\frac{\ell(\ell+2s-3)}{2}};q^{\ell(\ell-1)}\Big)_\infty\Big(-q^{\frac{\ell(\ell-2s+1)}{2}};q^{\ell(\ell-1)}\Big)_\infty\\
+q^{\frac{\ell^2-1}{8}}\,\frac{f_{2\ell(\ell-1)}^2}{f_{\ell(\ell-1)}^2}\Bigg].
\end{multline*}
\end{lemma}
\begin{proof}
See \cite[Theorem 1.2]{TBM}.
\end{proof}

\begin{lemma}\label{lem:2.10}
We have
\begin{align}
\frac{1}{f_1^4}&=\frac{f_4^{14}}{f_2^{14}f_8^4}+4q\,\frac{f_4^2f_8^4}{f_2^{10}},\label{eq:2.6}\\
\frac{f_2^2}{f_1}&=\frac{f_6f_9^2}{f_3f_{18}}+q\,\frac{f_{18}^2}{f_9},\label{eq:2.7}\\
\frac{f_2}{f_1^2}&=\frac{f_6^4f_9^6}{f_3^8f_{18}^3}+2q\,\frac{f_6^3f_9^3}{f_3^7}+4q^2\,\frac{f_6^2f_{18}^3}{f_3^6}.\label{eq:2.8}
\end{align}
\end{lemma}
\begin{proof}
See \cite{Hir}.
\end{proof}

\section{Proofs}

\begin{proposition}\label{prop:1}
For a positive integer $\ell\ge2$, and
\[
R_\ell(m)=\sum_{\substack{1\le j<m\\ \ell\nmid j}}j,
\]
we have
\begin{equation}\label{eq:3.1}
\sum_{\substack{m\ge1\\ \ell\nmid m}}\sum_{n=0}^\infty p_\ell^{\overline{\mathrm{mex}}}(m,n)z^mq^n=\frac{f_2f_\ell^2}{f_1^2f_{2\ell}}\sum_{\substack{m\ge1\\ \ell\nmid m}}z^mq^{R_\ell(m)}\big(1-q^m\big).
\end{equation}
\end{proposition}

\begin{proof}
Suppose $\omex_\ell(\pi)=m$, $\ell\nmid m$. Then
\begin{enumerate}[(1)]
\item every admissible integer $1\le r<m$, $\ell\nmid r$, must occur at least once as a non-overlined part;
\item the part $m$ is absent from the non-overlined parts;
\item every admissible integer greater than $m$ may occur as a non-overlined part with arbitrary multiplicity;
\item the overlined parts form an arbitrary partition into distinct admissible parts.
\end{enumerate}
The compulsory contribution of the non-overlined parts below $m$ is $q^{R_\ell(m)}\prod_{r<m,\ \ell\nmid r}\frac{1}{1-q^r}$. The unrestricted larger admissible non-overlined parts contribute $\prod_{r>m,\ \ell\nmid r}\frac{1}{1-q^r}$, and the overlined parts contribute $\prod_{r\ge1,\ \ell\nmid r}(1+q^r)$. The only admissible $r$ missing from the denominators is $r=m$. Thus
\begin{equation}\label{eq:3.2}
\sum_{n=0}^\infty p_\ell^{\overline{\mathrm{mex}}}(m,n)q^n=q^{R_\ell(m)}\big(1-q^m\big)\prod_{\substack{r\ge1\\ \ell\nmid r}}\frac{1}{1-q^r}\prod_{\substack{r\ge1\\ \ell\nmid r}}(1+q^r),\qquad\ell\nmid m.
\end{equation}
By the generating function of $\op_\ell(n)$ given in Section \ref{sec:intro},
\begin{equation}\label{eq:3.3}
\sum_{n=0}^\infty p_\ell^{\overline{\mathrm{mex}}}(m,n)q^n=q^{R_\ell(m)}\,\frac{f_2f_\ell^2}{f_1^2f_{2\ell}}\big(1-q^m\big),\qquad\ell\nmid m.
\end{equation}
Consequently, multiplying \eqref{eq:3.3} by $z^m$ and summing over $m\ge1$, $\ell\nmid m$, we have \eqref{eq:3.1}.
\end{proof}

\subsection{Proof of Theorem \ref{thm:1.1} and \ref{thm:1.2}.}

\begin{proof}[Proof of Theorem \ref{thm:1.1}]
Let $p_2^{\overline{\mathrm{mex}}}(r,n)$ be the number of $2$-regular overpartitions $\pi$ of $n$ such that $\omex_2(\pi)=2r-1$. We have
\begin{equation}\label{eq:3.4}
H_2(z,q):=\sum_{n\ge0}\sum_{r\ge1}p_2^{\overline{\mathrm{mex}}}(r,n)z^{2r-1}q^n=\sum_{r=1}^\infty z^{2r-1}\sum_{n=0}^\infty p_2^{\overline{\mathrm{mex}}}(r,n)q^n.
\end{equation}
Since $\omex_2(\pi)=2r-1$, the compulsory non-overlined parts are $1,3,5,\ldots,2r-3$, and their sum is
\[
1+3+\cdots+(2r-3)=(r-1)^2.
\]
Hence these compulsory parts contribute $q^{(r-1)^2}\prod_{j=1}^{r-1}\frac{1}{1-q^{2j-1}}$. The non-overlined part $2r-1$ is excluded, while every larger odd part may occur as a non-overlined part with arbitrary multiplicity. The overlined parts form a partition into distinct odd parts, which contributes $(-q;q^2)_\infty$. Thus,
\begin{equation}\label{eq:3.5}
\sum_{n=0}^\infty p_2^{\overline{\mathrm{mex}}}(r,n)q^n=q^{(r-1)^2}\prod_{j=1}^{r-1}\frac{1}{1-q^{2j-1}}\prod_{j=r+1}^{\infty}\frac{1}{1-q^{2j-1}}\,(-q;q^2)_\infty=\frac{q^{(r-1)^2}\big(1-q^{2r-1}\big)(-q;q^2)_\infty}{(q;q^2)_\infty}.
\end{equation}
Combining \eqref{eq:3.4} and \eqref{eq:3.5}, we get
\begin{equation}\label{eq:3.6}
H_2(z,q)=\frac{(-q;q^2)_\infty}{(q;q^2)_\infty}\sum_{r=1}^\infty z^{2r-1}q^{(r-1)^2}\big(1-q^{2r-1}\big).
\end{equation}
Differentiating $H_2(z,q)$ with respect to $z$ and putting $z=1$, we get
\begin{align}
\sum_{n\ge0}\sm{2}(n)q^n&=\frac{\partial}{\partial z}H_2(z,q)\Big|_{z=1}\notag\\
&=\frac{(-q;q^2)_\infty}{(q;q^2)_\infty}\sum_{r=1}^\infty(2r-1)q^{(r-1)^2}\big(1-q^{2r-1}\big)\notag\\
&=\frac{(-q;q^2)_\infty}{(q;q^2)_\infty}\left[\sum_{r=1}^\infty(2r-1)q^{(r-1)^2}-\sum_{r=1}^\infty(2r-1)q^{(r-1)^2}q^{2r-1}\right]\notag\\
&=\frac{(-q;q^2)_\infty}{(q;q^2)_\infty}\left[\sum_{t=0}^\infty(2t+1)q^{t^2}-\sum_{t=0}^\infty(2t+1)q^{(t+1)^2}\right]\notag\\
&=\frac{(-q;q^2)_\infty}{(q;q^2)_\infty}\left[\sum_{t=0}^\infty(2t+1)q^{t^2}-\sum_{t=1}^\infty(2t-1)q^{t^2}\right]\notag\\
&=\frac{(-q;q^2)_\infty}{(q;q^2)_\infty}\left[1+2\sum_{t=1}^\infty q^{t^2}\right]=\frac{(-q;q^2)_\infty}{(q;q^2)_\infty}\sum_{t\in\mathbb{Z}}q^{t^2}.\label{eq:3.7}
\end{align}
Since $(-q;q^2)_\infty=\frac{f_2^2}{f_1f_4}$ and $(q;q^2)_\infty=\frac{f_1}{f_2}$, Lemma \ref{lem:2.1} gives
\begin{equation}\label{eq:3.8}
\sum_{n\ge0}\sm{2}(n)q^n=\frac{f_2^3}{f_1^2f_4}\cdot\frac{f_2^5}{f_1^2f_4^2}=\frac{f_2^8}{f_1^4f_4^3}.
\end{equation}
\end{proof}

In order to prove Theorem \ref{thm:1.2}, we need the following proposition.

\begin{proof}[Proof of Theorem \ref{thm:1.2}]
From Proposition \ref{prop:1}, we have
\begin{equation}\label{eq:3.9}
H_3(z,q):=\sum_{n\ge0}\sum_{\substack{m\ge1\\ 3\nmid m}}p_3^{\overline{\mathrm{mex}}}(m,n)z^mq^n=\frac{f_2f_3^2}{f_1^2f_6}\sum_{\substack{m\ge1\\ 3\nmid m}}z^mq^{R_3(m)}\big(1-q^m\big).
\end{equation}
Using Lemma \ref{lem:2.8} in \eqref{eq:3.9}, we get
\begin{equation}\label{eq:3.10}
H_3(z,q)=\frac{f_2f_3^2}{f_1^2f_6}\left[\sum_{r\ge0}z^{3r+1}q^{3r^2}\big(1-q^{3r+1}\big)+\sum_{r\ge0}z^{3r+2}q^{3r^2+3r+1}\big(1-q^{3r+2}\big)\right].
\end{equation}
Differentiating \eqref{eq:3.10} with respect to $z$ and putting $z=1$, we obtain
\begin{equation}\label{eq:3.11}
\sum_{n\ge0}\sm{3}(n)q^n=\frac{\partial}{\partial z}H_3(z,q)\Big|_{z=1}=\frac{f_2f_3^2}{f_1^2f_6}\,N_3(q),
\end{equation}
where
\begin{align}
N_3(q)&=\sum_{r\ge0}(3r+1)q^{3r^2}\big(1-q^{3r+1}\big)+\sum_{r\ge0}(3r+2)q^{3r^2+3r+1}\big(1-q^{3r+2}\big)\label{eq:3.12}\\
&=\sum_{r\ge0}(3r+1)q^{3r^2}+\sum_{r\ge0}q^{3r^2+3r+1}-\sum_{r\ge0}(3r+2)q^{3(r+1)^2}\notag\\
&=\sum_{r\ge0}(3r+1)q^{3r^2}-\sum_{r\ge1}(3r-1)q^{3r^2}+\sum_{r\ge0}q^{3r^2+3r+1}\notag\\
&=1+2\sum_{r\ge1}q^{3r^2}+\sum_{r\ge0}q^{3r^2+3r+1}=\sum_{r\in\mathbb{Z}}q^{3r^2}+\sum_{r\ge0}q^{3r^2+3r+1}.\label{eq:3.13}
\end{align}
Here, in the second line we combined $-(3r+1)q^{3r^2+3r+1}+(3r+2)q^{3r^2+3r+1}=q^{3r^2+3r+1}$, and used $3r^2+3r+1+(3r+2)=3(r+1)^2$.

Note that, for any real number $r$,
\begin{equation}\label{eq:3.14}
(-r-1)^2+(-r-1)=r^2+r.
\end{equation}
Using \eqref{eq:3.14}, we get
\begin{equation}\label{eq:3.15}
\sum_{r\ge0}q^{3r^2+3r+1}=\frac{1}{2}\sum_{r\in\mathbb{Z}}q^{3r^2+3r+1}=\frac{q}{2}\sum_{r=-\infty}^{\infty}\big(q^3\big)^{r^2}\big(q^3\big)^{r}.
\end{equation}
Applying $Q=q^3$ and $z=q^3$ in Lemma \ref{lem:2.2}, we have
\begin{equation}\label{eq:3.16}
\sum_{r=-\infty}^{\infty}q^{3r^2+3r}=\big(-q^6;q^6\big)_\infty\big(-1;q^6\big)_\infty\big(q^6;q^6\big)_\infty=2\big(-q^6;q^6\big)_\infty^2\big(q^6;q^6\big)_\infty,
\end{equation}
since $(-1;q^6)_\infty=2(-q^6;q^6)_\infty$. Further observe that
\begin{equation}\label{eq:3.17}
\big(-q^6;q^6\big)_\infty=\frac{(q^{12};q^{12})_\infty}{(q^6;q^6)_\infty}.
\end{equation}
Using \eqref{eq:3.16} and \eqref{eq:3.17} in \eqref{eq:3.15}, we get
\begin{equation}\label{eq:3.18}
\sum_{r\ge0}q^{3r^2+3r+1}=q\,\frac{(q^{12};q^{12})_\infty^2}{(q^6;q^6)_\infty}.
\end{equation}
Also, using Lemma \ref{lem:2.1} with $q\to q^3$, we have
\begin{equation}\label{eq:3.19}
\sum_{r\in\mathbb{Z}}q^{3r^2}=\frac{(q^6;q^6)_\infty^5}{(q^3;q^3)_\infty^2(q^{12};q^{12})_\infty^2}.
\end{equation}
Thus, by \eqref{eq:3.11}, \eqref{eq:3.13}, \eqref{eq:3.18} and \eqref{eq:3.19},
\begin{align}
\sum_{n\ge0}\sm{3}(n)q^n&=\frac{f_2f_3^2}{f_1^2f_6}\left[\frac{f_6^5}{f_3^2f_{12}^2}+q\frac{f_{12}^2}{f_6}\right]\label{eq:3.20}\\
&=\frac{f_2f_6^4}{f_1^2f_{12}^2}+q\,\frac{f_2f_3^2f_{12}^2}{f_1^2f_6^2}.\notag
\end{align}
\end{proof}

\subsection{Proof of Theorem \ref{thm:1.3}.}

\begin{proof}
Differentiating \eqref{eq:3.1} with respect to $z$ and putting $z=1$, we get
\begin{equation}\label{eq:3.21}
\sum_{n\ge0}\sm{\ell}(n)q^n=\frac{f_2f_\ell^2}{f_1^2f_{2\ell}}\sum_{\substack{m\ge1\\ \ell\nmid m}}m\,q^{R_\ell(m)}\big(1-q^m\big).
\end{equation}
For $\ell$-regular partitions there are no overlined parts, so the argument in the proof of Proposition \ref{prop:1}, with the factor
\[
\prod_{r\ge1,\ \ell\nmid r}(1+q^r)
\]
is removed. Therefore
\[
\sum_{n\ge0}\sigma_\ell^{\mathrm{mex}}(n)q^n=\frac{f_\ell}{f_1}\sum_{\substack{m\ge1\\ \ell\nmid m}}m\,q^{R_\ell(m)}\big(1-q^m\big).
\]
Comparing with \eqref{eq:3.21}, we obtain
\begin{equation}\label{eq:3.22}
\sum_{n\ge0}\sm{\ell}(n)q^n=\frac{f_2f_\ell}{f_1f_{2\ell}}\sum_{n\ge0}\sigma_\ell^{\mathrm{mex}}(n)q^n.
\end{equation}
Substituting Lemma \ref{lem:2.9} in \eqref{eq:3.22}, and using
\[
\frac{f_2f_\ell}{f_1f_{2\ell}}\cdot\frac{f_\ell\,f_{\ell(\ell-1)}}{f_1}=\frac{f_2f_\ell^2f_{\ell(\ell-1)}}{f_1^2f_{2\ell}},
\]
we obtain \eqref{eq:1.4} for even $\ell$ and \eqref{eq:1.5} for odd $\ell$.
\end{proof}

\begin{proof}[Proof of Theorem \ref{thm:1.4}]
Consider the map $\phi:\cP_\ell(n)\to\cP_\ell(n+1)$ defined by $\phi(\pi)=\pi\cup\{1\}$, where the new part $1$ is non-overlined. Since $1$ is never divisible by $\ell$, $\phi(\pi)$ is still an $\ell$-regular overpartition, and $|\phi(\pi)|=|\pi|+1=n+1$. If $\pi$ already contains $\overline{1}$, the new part is simply an additional non-overlined $1$, which is allowed. The map is injective: if $\phi(\pi_1)=\phi(\pi_2)$, then removing one non-overlined $1$ from both gives $\pi_1=\pi_2$.

Now we consider two cases. Firstly, suppose that $1$ is a non-overlined part of $\pi$. Then $\omex_\ell(\pi)>1$. Adding another non-overlined $1$ does not change whether any other part occurs as a non-overlined part. Hence $\omex_\ell(\phi(\pi))=\omex_\ell(\pi)$. Secondly, suppose that $1$ is not a non-overlined part of $\pi$. Then $\omex_\ell(\pi)=1$, and after adding $1$, $1$ is a non-overlined part of $\phi(\pi)$. Therefore $\omex_\ell(\phi(\pi))$ is the smallest positive integer other than $1$ that is not divisible by $\ell$ and does not occur as a non-overlined part. Thus $\omex_\ell(\phi(\pi))>1=\omex_\ell(\pi)$. In both cases $\omex_\ell(\phi(\pi))\ge\omex_\ell(\pi)$. Since $\phi$ is injective, we have
\[
\sm{\ell}(n+1)=\sum_{\mu\in\cP_\ell(n+1)}\omex_\ell(\mu)\ge\sum_{\pi\in\cP_\ell(n)}\omex_\ell(\phi(\pi))\ge\sum_{\pi\in\cP_\ell(n)}\omex_\ell(\pi)=\sm{\ell}(n).
\]
Hence $\sm{\ell}(n+1)\ge\sm{\ell}(n)$ for all $n\ge0$. Therefore, the sequence $\{\sm{\ell}(n)\}_{n\ge0}$ is weakly increasing for every $\ell\ge2$.

Every overpartition $\pi\in\cP_\ell(n)$ without a non-overlined $1$ satisfies $\omex_\ell(\pi)=1$. After adding a part $1$, its mex becomes strictly larger, i.e. $\omex_\ell(\phi(\pi))>1=\omex_\ell(\pi)$. If $1$ is a non-overlined part of $\pi$, then $\omex_\ell(\phi(\pi))=\omex_\ell(\pi)$. Hence
\begin{equation}\label{eq:3.23}
\sm{\ell}(n+1)-\sm{\ell}(n)\ge\#\{\pi\in\cP_\ell(n):\ 1\text{ is not a non-overlined part of }\pi\}.
\end{equation}
Now consider the overpartitions in $\cP_\ell(n)$ that contain at least one non-overlined part $1$. Removing one such $1$ gives a bijection with $\cP_\ell(n-1)$. Therefore,
\begin{equation}\label{eq:3.24}
\#\{\pi\in\cP_\ell(n):\ 1\text{ is not a non-overlined part of }\pi\}=\op_\ell(n)-\op_\ell(n-1).
\end{equation}
Comparing \eqref{eq:3.23} and \eqref{eq:3.24}, we get $\sm{\ell}(n+1)-\sm{\ell}(n)\ge\op_\ell(n)-\op_\ell(n-1)$.
\end{proof}

\begin{proof}[Proof of Theorem \ref{thm:1.5}]
Let
\[
G(q)=\sum_{n\ge0}\sm{2}(n)q^n=\frac{(q^2;q^2)_\infty^8}{(q;q)_\infty^4(q^4;q^4)_\infty^3}.
\]
Applying the transformation formula for Dedekind's eta function \cite{Kob}, we have
\begin{equation}\label{eq:3.25}
(e^{-y};e^{-y})_\infty\sim\sqrt{\frac{2\pi}{y}}\exp\Big(-\frac{\pi^2}{6y}\Big),\qquad y\to0^+.
\end{equation}
Thus
\[
(e^{-2y};e^{-2y})_\infty^8\sim\Big(\frac{\pi}{y}\Big)^4e^{-\frac{2\pi^2}{3y}},\qquad(e^{-y};e^{-y})_\infty^4\sim\Big(\frac{2\pi}{y}\Big)^2e^{-\frac{2\pi^2}{3y}},
\]
and
\[
(e^{-4y};e^{-4y})_\infty^3\sim\Big(\frac{\pi}{2y}\Big)^{3/2}e^{-\frac{\pi^2}{8y}}.
\]
Therefore,
\[
G(e^{-y})\sim\frac{\big(\frac{\pi}{y}\big)^4}{\big(\frac{2\pi}{y}\big)^2\big(\frac{\pi}{2y}\big)^{3/2}}\,e^{-\frac{2\pi^2}{3y}+\frac{2\pi^2}{3y}+\frac{\pi^2}{8y}}=\sqrt{\frac{\pi}{2}}\,y^{-1/2}e^{\frac{\pi^2}{8y}}.
\]
From Theorem \ref{thm:1.4}, we know that the sequence $\{\sm{2}(n)\}$ is weakly increasing. Thus, by applying Lemma \ref{lem:2.3} with parameters $\beta=\sqrt{\frac{\pi}{2}}$, $\gamma=-\frac{1}{2}$, $\mu=\frac{\pi^2}{8}$, we obtain
\[
\sm{2}(n)\sim\frac{\sqrt{\pi/2}}{2\sqrt{\pi n}}\,e^{2\sqrt{\frac{\pi^2n}{8}}}=\frac{1}{2\sqrt{2n}}\,e^{\pi\sqrt{n/2}},\qquad\text{as } n\to\infty.\qedhere
\]
\end{proof}

\begin{proof}[Proof of Theorem \ref{thm:1.6}]
From \eqref{eq:3.20}, we have
\[
\sum_{n\ge0}\sm{3}(n)q^n=\frac{(q^2;q^2)_\infty(q^6;q^6)_\infty^4}{(q;q)_\infty^2(q^{12};q^{12})_\infty^2}+q\,\frac{(q^2;q^2)_\infty(q^3;q^3)_\infty^2(q^{12};q^{12})_\infty^2}{(q;q)_\infty^2(q^6;q^6)_\infty^2}.
\]
Let
\[
S_1(q):=\frac{(q^2;q^2)_\infty(q^6;q^6)_\infty^4}{(q;q)_\infty^2(q^{12};q^{12})_\infty^2}.
\]
Using \eqref{eq:3.25}, we get
\[
(e^{-2y};e^{-2y})_\infty(e^{-6y};e^{-6y})_\infty^4\sim\sqrt{\frac{\pi}{y}}\Big(\frac{\pi}{3y}\Big)^2e^{-\frac{\pi^2}{12y}-\frac{\pi^2}{9y}},
\]
and
\[
(e^{-y};e^{-y})_\infty^2(e^{-12y};e^{-12y})_\infty^2\sim\frac{2\pi}{y}\cdot\frac{\pi}{6y}\,e^{-\frac{\pi^2}{3y}-\frac{\pi^2}{36y}}.
\]
Since $-\frac{1}{12}-\frac{1}{9}+\frac{1}{3}+\frac{1}{36}=\frac{1}{6}$, we obtain
\[
S_1(e^{-y})\sim\frac{\sqrt{\pi}}{3}\,y^{-1/2}e^{\frac{\pi^2}{6y}}.
\]
Similarly, let
\[
S_2(q):=q\,\frac{(q^2;q^2)_\infty(q^3;q^3)_\infty^2(q^{12};q^{12})_\infty^2}{(q;q)_\infty^2(q^6;q^6)_\infty^2}.
\]
Using \eqref{eq:3.25}, we get
\[
(e^{-2y};e^{-2y})_\infty(e^{-3y};e^{-3y})_\infty^2(e^{-12y};e^{-12y})_\infty^2\sim\sqrt{\frac{\pi}{y}}\cdot\frac{2\pi}{3y}\cdot\frac{\pi}{6y}\,e^{-\frac{\pi^2}{12y}-\frac{\pi^2}{9y}-\frac{\pi^2}{36y}},
\]
and
\[
(e^{-y};e^{-y})_\infty^2(e^{-6y};e^{-6y})_\infty^2\sim\frac{2\pi}{y}\cdot\frac{\pi}{3y}\,e^{-\frac{\pi^2}{3y}-\frac{\pi^2}{18y}}.
\]
Since $-\frac{1}{12}-\frac{1}{9}-\frac{1}{36}+\frac{1}{3}+\frac{1}{18}=\frac{1}{6}$, and since $q\to1$ as $y\to0^+$, we get
\[
S_2(e^{-y})\sim\frac{\sqrt{\pi}}{6}\,y^{-1/2}e^{\frac{\pi^2}{6y}}.
\]
Note that both terms have the same exponential growth. Consequently,
\[
S_1(e^{-y})+S_2(e^{-y})=\sum_{n\ge0}\sm{3}(n)e^{-ny}\sim\Big(\frac{\sqrt{\pi}}{3}+\frac{\sqrt{\pi}}{6}\Big)y^{-1/2}e^{\frac{\pi^2}{6y}}=\frac{\sqrt{\pi}}{2}\,y^{-1/2}e^{\frac{\pi^2}{6y}}.
\]
Applying Lemma \ref{lem:2.3} with parameters $\beta=\frac{\sqrt{\pi}}{2}$, $\gamma=-\frac{1}{2}$, $\mu=\frac{\pi^2}{6}$, we get
\[
\sm{3}(n)\sim\frac{\sqrt{\pi}/2}{2\sqrt{\pi n}}\,e^{2\sqrt{\frac{\pi^2n}{6}}}=\frac{1}{4\sqrt{n}}\,e^{\pi\sqrt{2n/3}},\qquad\text{as } n\to\infty.\qedhere
\]
\end{proof}

\subsection{Proof of Theorems \ref{thm:1.7}, \ref{thm:1.8} and \ref{thm:parity-general}.}
In order to prove Theorem \ref{thm:1.7}, we need the following result.

\begin{proposition}\label{prop:2}
For every $\ell\ge2$, we have
\begin{equation}\label{eq:3.26}
\overline{P}_\ell(q)=\frac{f_2f_\ell^2}{f_1^2f_{2\ell}}\equiv1\pmod2.
\end{equation}
\end{proposition}
\begin{proof}
By Lemma \ref{lem:2.4} with $p=2$ and $j=1$, $f_1^2\equiv f_2\pmod 2$ and $f_\ell^2\equiv f_{2\ell}\pmod2$. Hence $\frac{f_2f_\ell^2}{f_1^2f_{2\ell}}\equiv\frac{f_2f_{2\ell}}{f_2f_{2\ell}}=1\pmod2$.
\end{proof}

\begin{proof}[Proof of Theorem \ref{thm:1.7}]
(1) From \eqref{eq:3.7} and \eqref{eq:3.8}, we have
\begin{equation}\label{eq:3.27}
\sum_{n=0}^\infty\sm{2}(n)q^n=\overline{P}_2(q)\sum_{t\in\mathbb{Z}}q^{t^2}=\overline{P}_2(q)\Big(1+2\sum_{t\ge1}q^{t^2}\Big).
\end{equation}
By Proposition \ref{prop:2},
\[
\sum_{n=0}^\infty\sm{2}(n)q^n\equiv1\pmod2.
\]
Hence $\sm{2}(n)$ is even for every $n\ge1$.

(2) From \eqref{eq:3.11}, \eqref{eq:3.13} and Proposition \ref{prop:2}, we have
\begin{equation}\label{eq:3.28}
\sum_{n=0}^\infty\sm{3}(n)q^n\equiv\sum_{r\in\mathbb{Z}}q^{3r^2}+\sum_{r\ge0}q^{3r^2+3r+1}\equiv1+\sum_{r\ge0}q^{3r^2+3r+1}\pmod2,
\end{equation}
since $\sum_{r\in\mathbb{Z}}q^{3r^2}=1+2\sum_{r\ge1}q^{3r^2}$. The exponents $3r^2+3r+1$ $(r\ge0)$ are distinct positive integers. Comparing the coefficients of $q^n$ on both sides of \eqref{eq:3.28}, we obtain
\[
\sm{3}(n)\equiv1\pmod2\iff n=0\ \text{ or }\ n=3k^2+3k+1\ \text{ for some integer } k\ge0.\qedhere
\]
\end{proof}

\begin{proof}[Proof of Theorem \ref{thm:1.8}]
For a real number $X\ge1$, we define
\[
M_1(X):=\#\{n\le X:\ \sm{3}(n)\text{ is odd}\},\qquad M_2(X):=\#\{n\le X:\ \sm{3}(n)\text{ is even}\}.
\]
Note that, for any real number $k\ge0$,
\[
3k^2+3k+1\le X\iff k\le\frac{-3+\sqrt{12X-3}}{6}.
\]
Using Theorem \ref{thm:1.7}(2), we get
\begin{equation}\label{eq:3.29}
M_1(X)=\left\lfloor\frac{-3+\sqrt{12X-3}}{6}\right\rfloor+1=\frac{\sqrt{12X-3}}{6}+O(1).
\end{equation}
For large $x$, we have $\sqrt{x+a}=\sqrt{x}+\frac{a}{2\sqrt{x}}+O\big(x^{-3/2}\big)$. In particular,
\[
\frac{\sqrt{12X-3}}{6}=\frac{\sqrt{12X}}{6}+O\big(X^{-1/2}\big)=\sqrt{\frac{X}{3}}+O\big(X^{-1/2}\big).
\]
Thus
\[
M_1(X)=\sqrt{\frac{X}{3}}+O(1).
\]
Therefore,
\[
M_2(X)=\lfloor X\rfloor-M_1(X)=X-\sqrt{\frac{X}{3}}+O(1).\qedhere
\]
\end{proof}

\begin{proof}[Proof of Theorem \ref{thm:parity-general}]
Throughout the proof, $\ell\ge2$ is fixed and we write $R(m)=R_\ell(m)$. Let
\[
1=m_1<m_2<m_3<\cdots
\]
be the positive integers not divisible by $\ell$, listed in increasing order.

By \eqref{eq:3.21},
\begin{equation}\label{eq:pg1}
\sum_{n\ge0}\sm{\ell}(n)q^n=\overline{P}_\ell(q)\,S_\ell(q),\qquad S_\ell(q):=\sum_{\substack{m\ge1\\ \ell\nmid m}}m\,q^{R(m)}\big(1-q^m\big)=\sum_{j\ge1}m_j\,q^{R(m_j)}\big(1-q^{m_j}\big).
\end{equation}
By Proposition \ref{prop:2}, $\overline{P}_\ell(q)\equiv1\pmod2$, and therefore
\begin{equation}\label{eq:pg2}
\sum_{n\ge0}\sm{\ell}(n)q^n\equiv S_\ell(q)\pmod2.
\end{equation}

Next, for every $j\ge1$, the integer $m_j$ is the only positive integer in the interval $[m_j,m_{j+1})$ that is not divisible by $\ell$. Hence, by the definition of $R$,
\begin{equation}\label{eq:pg3}
R(m_{j+1})=\sum_{\substack{1\le i<m_{j+1}\\ \ell\nmid i}}i=\sum_{\substack{1\le i<m_j\\ \ell\nmid i}}i+m_j=R(m_j)+m_j\qquad(j\ge1).
\end{equation}
Moreover, $R(m_1)=R(1)=0$, and the sequence $\{R(m_j)\}_{j\ge1}$ is strictly increasing. Using \eqref{eq:pg3} in \eqref{eq:pg1}, we obtain
\begin{align*}
S_\ell(q)&=\sum_{j\ge1}m_j\,q^{R(m_j)}-\sum_{j\ge1}m_j\,q^{R(m_j)+m_j}
=\sum_{j\ge1}m_j\,q^{R(m_j)}-\sum_{j\ge1}m_j\,q^{R(m_{j+1})}\\
&=\sum_{j\ge1}m_j\,q^{R(m_j)}-\sum_{j\ge2}m_{j-1}\,q^{R(m_j)}
=m_1q^{R(m_1)}+\sum_{j\ge2}\big(m_j-m_{j-1}\big)q^{R(m_j)}.
\end{align*}
This gives
\begin{equation}\label{eq:pg4}
S_\ell(q)=1+\sum_{j\ge2}\big(m_j-m_{j-1}\big)q^{R(m_j)}.
\end{equation}

We now determine the differences $m_j-m_{j-1}$. Let $j\ge2$. If $m_j\not\equiv1\pmod\ell$, then $\ell\nmid m_j-1$, so $m_{j-1}=m_j-1$ and $m_j-m_{j-1}=1$. If $m_j\equiv1\pmod\ell$, then $\ell\mid m_j-1$, while $m_j-2\equiv-1\not\equiv0\pmod\ell$ since $\ell\ge2$. Hence $m_{j-1}=m_j-2$ and $m_j-m_{j-1}=2$. Therefore
\[
m_j-m_{j-1}\equiv\begin{cases}1\pmod2,& m_j\not\equiv1\pmod\ell,\\ 0\pmod2,& m_j\equiv1\pmod\ell.\end{cases}
\]
Combining this with \eqref{eq:pg2} and \eqref{eq:pg4}, we get
\begin{equation}\label{eq:pg5}
\sum_{n\ge0}\sm{\ell}(n)q^n\equiv1+\sum_{\substack{j\ge2\\ m_j\not\equiv1\ (\mathrm{mod}\ \ell)}}q^{R(m_j)}\pmod2.
\end{equation}
Since the exponents $R(m_j)$, $j\ge1$, are pairwise distinct and $R(m_j)>0$ for $j\ge2$, no two terms on the right-hand side of \eqref{eq:pg5} have the same power of $q$. Comparing the coefficients of $q^n$ on both sides of \eqref{eq:pg5}, we conclude that $\sm{\ell}(n)$ is odd if and only if $n=0$ or $n=R(m_j)$ for some $j\ge2$ with $m_j\not\equiv1\pmod\ell$.

 Note that the integers $m_j$ with $j\ge2$ and $m_j\not\equiv1\pmod\ell$ are precisely the integers $\ell k+s$ with $k\ge0$ and $2\le s\le\ell-1$. By Lemma \ref{lem:2.8},
\[
R(\ell k+s)=\frac{\ell(\ell-1)}{2}k^2+(s-1)\ell k+\frac{s(s-1)}{2}.
\]
This proves the first assertion. 

It remains to count the odd values. Fixing $s\in\{2,\ldots,\ell-1\}$ and putting
\[
a=\frac{\ell(\ell-1)}{2},\qquad b=(s-1)\ell,\qquad c=\frac{s(s-1)}{2}, 
\]
we see that the map $k\mapsto R(\ell k+s)=ak^2+bk+c$ is strictly increasing on $k\ge0$. For $X\ge c$, the inequality $ak^2+bk+c\le X$ with $k\ge0$ is equivalent to $0\le k\le\kappa_s(X)$, where
\[
\kappa_s(X)=\frac{-b+\sqrt{b^2-4ac+4aX}}{2a}.
\]
Hence the number $N_s(X)$ of integers $k\ge0$ with $R(\ell k+s)\le X$ equals $\lfloor\kappa_s(X)\rfloor+1$ for $X\ge c$, and $N_s(X)=0$ for $X<c$. Since $b^2-4ac$ does not depend on $X$, we have
\[
\sqrt{b^2-4ac+4aX}=\sqrt{4aX}+O\big(X^{-1/2}\big),
\]
and therefore
\[
N_s(X)=\frac{\sqrt{4aX}}{2a}+O(1)=\sqrt{\frac{X}{a}}+O(1)=\sqrt{\frac{2X}{\ell(\ell-1)}}+O(1).
\]
By the first part of the theorem, the values $R(\ell k+s)$ for different pairs $(k,s)$ are distinct, and they are distinct from $0$. Hence
\[
\#\{0\le n\le X:\ \sm{\ell}(n)\ \text{is odd}\}=1+\sum_{s=2}^{\ell-1}N_s(X)=(\ell-2)\sqrt{\frac{2X}{\ell(\ell-1)}}+O(1).
\]
For $\ell\ge3$ the main term is a positive multiple of $\sqrt{X}$, so $\sm{\ell}(n)$ is odd for infinitely many $n$, and the number of such $n\le X$ is $\gg\sqrt{X}$. This completes the proof.
\end{proof}

\subsection{Proof of Theorem \ref{thm:1.9}.}

\begin{proof}
(1) By Lemma \ref{lem:2.4} with $p=3$ and $j=1$, we have $f_1^3\equiv f_3$, $f_2^6\equiv f_6^2$ and $f_4^3\equiv f_{12}\pmod3$. Hence, by Theorem \ref{thm:1.1},
\[
\sum_{n\ge0}\sm{2}(n)q^n=\frac{f_2^8}{f_1^4f_4^3}\equiv\frac{f_6^2}{f_3f_{12}}\cdot\frac{f_2^2}{f_1}\pmod3.
\]
Using \eqref{eq:2.7}, we get
\begin{equation}\label{eq:3.30}
\sum_{n\ge0}\sm{2}(n)q^n\equiv\frac{f_6^3f_9^2}{f_3^2f_{12}f_{18}}+q\,\frac{f_6^2f_{18}^2}{f_3f_9f_{12}}\pmod3.
\end{equation}
The right-hand side of \eqref{eq:3.30} contains no term $q^{3n+2}$. Therefore $\sm{2}(3n+2)\equiv0\pmod3$.

(2) Substituting \eqref{eq:2.6} in Theorem \ref{thm:1.1}, we obtain
\begin{equation}\label{eq:3.31}
\sum_{n\ge0}\sm{2}(n)q^n=\frac{f_4^{11}}{f_2^6f_8^4}+4q\,\frac{f_8^4}{f_2^2f_4}.
\end{equation}
Extracting the terms with even powers of $q$ and replacing $q^2$ by $q$, we get
\[
\sum_{n\ge0}\sm{2}(2n)q^n=\frac{f_2^{11}}{f_1^6f_4^4}=\frac{f_2^5}{f_1^2f_4^2}\cdot\frac{f_2^6}{f_1^4f_4^2}.
\]
By Lemma \ref{lem:2.4} with $p=2$ and $j=2$, $f_1^4\equiv f_2^2$ and $f_2^4\equiv f_4^2\pmod4$, so $\frac{f_2^6}{f_1^4f_4^2}\equiv\frac{f_2^4}{f_4^2}\equiv1\pmod4$. Hence, by Lemma \ref{lem:2.1},
\begin{equation}\label{eq:3.32}
\sum_{n\ge0}\sm{2}(2n)q^n\equiv\frac{f_2^5}{f_1^2f_4^2}=\sum_{k\in\mathbb{Z}}q^{k^2}=1+2\sum_{k\ge1}q^{k^2}\pmod4.
\end{equation}
Also, by \eqref{eq:3.31}, $\sm{2}(2n+1)\equiv0\pmod4$. Therefore, for $n\ge1$, $\sm{2}(n)\equiv2\pmod4$ if $n=2k^2$ for some $k\ge1$, and $\sm{2}(n)\equiv0\pmod4$ otherwise.

(3) Extracting the terms with odd powers of $q$ in \eqref{eq:3.31}, dividing by $q$ and replacing $q^2$ by $q$, we get
\begin{equation}\label{eq:3.33}
\sum_{n\ge0}\sm{2}(2n+1)q^n=4\frac{f_4^4}{f_1^2f_2}.
\end{equation}
By Lemma \ref{lem:2.4}, $f_1^2\equiv f_2$ and $f_2^2\equiv f_4\pmod2$, so by Lemma \ref{lem:2.6},
\[
\frac{f_4^4}{f_1^2f_2}\equiv f_4^3=\sum_{k\ge0}(-1)^k(2k+1)q^{2k(k+1)}\equiv\sum_{k\ge0}q^{2k(k+1)}\pmod2.
\]
Hence $\sum_n\sm{2}(2n+1)q^n\equiv4\sum_{k\ge0}q^{2k(k+1)}\pmod8$. Since $n=2k(k+1)$ if and only if $2n+1=(2k+1)^2$, the claim follows.

(4) By Theorem \ref{thm:1.2},
\[
\sum_{n\ge0}\sm{3}(n)q^n=\frac{f_6^4}{f_{12}^2}\cdot\frac{f_2}{f_1^2}+q\,\frac{f_3^2f_{12}^2}{f_6^2}\cdot\frac{f_2}{f_1^2}.
\]
Substituting \eqref{eq:2.8}, extracting the terms with powers $q^{3n}$ and replacing $q^3$ by $q$, we get
\[
\sum_{n\ge0}\sm{3}(3n)q^n=\frac{f_2^8f_3^6}{f_1^8f_4^2f_6^3}+4q\,\frac{f_4^2f_6^3}{f_1^4}.
\]
By Lemma \ref{lem:2.4} with $p=2$ and $j=2$, $f_1^8\equiv f_2^4$, $f_2^4\equiv f_4^2$ and $f_3^4\equiv f_6^2\pmod4$. Hence, by Lemma \ref{lem:2.1} with $q$ replaced by $-q^3$ (note that $(-q;-q)_\infty=\frac{f_2^3}{f_1f_4}$),
\[
\sum_{n\ge0}\sm{3}(3n)q^n\equiv\frac{f_3^2}{f_6}=\sum_{k\in\mathbb{Z}}(-1)^kq^{3k^2}\equiv1+2\sum_{k\ge1}q^{3k^2}\pmod4.
\]
Therefore, for $n\ge1$, $\sm{3}(3n)\equiv2\pmod4$ if $n=3k^2$ for some $k\ge1$, and $\sm{3}(3n)\equiv0\pmod4$ otherwise. Since $3n+1$ and $3n+2$ are never of the form $3k^2$, we get $\sm{3}(9n+3)\equiv\sm{3}(9n+6)\equiv0\pmod4$.
\end{proof}

\noindent\textbf{Data availability statement:} There is no data associated with our manuscript.

\end{document}